\documentclass[leqno]{article}  
\usepackage{epic,eepic,amsmath,amstext,amsthm,amssymb}   
\numberwithin{equation}{section}

\allowdisplaybreaks
\swapnumbers
\theoremstyle{plain}

\def\seq{\lesssim }

\def\bthm#1.#2 #3\ethm{
\begin{\ifx#1ttheorem\fi%
\ifx#1llemma\fi% 
\ifx#1ccorollary\fi% 
\ifx#1pproposition\fi%
\ifx#1ddefinition\fi}
\label{#1.#2}  
#3 \end{\ifx#1ttheorem\fi%
\ifx#1llemma\fi%
\ifx#1ccorollary\fi% 
\ifx#1pproposition\fi%
\ifx#1ddefinition\fi}}

\def\t#1/{theorem~\ref{t#1}}   \def\T#1/{Theorem~\ref{t#1}} 
\def\c#1/{corollary~\ref{c#1}}   \def\C#1/{Corollary~\ref{c#1}} 
\def\l#1/{lemma~\ref{l#1}}        \def\L#1/{Lemma~\ref{l#1}}  
\def\s#1/{section~\ref{s#1}}      
\def\e#1/{(\ref{e#1})}
\def\d#1/{definition~\ref{d#1}}
\def\f#1/{figure~\ref{f#1}}
\def\Label #1 {\label{#1}}

\def\norm#1.#2.{\lVert#1\rVert_{#2}}
\def\Norm#1.#2.{\bigl\lVert#1\bigr\rVert_{#2}}
\def\NOrm#1.#2.{\Bigl\lVert#1\Bigr\rVert_{#2}}
\def\NORm#1.#2.{\biggl\lVert#1\biggr\rVert_{#2}}
\def\NORM#1.#2.{\Biggl\lVert#1\Biggr\rVert_{#2}}

\def\ip#1,#2.{\langle #1,#2\rangle}
\def\Ip#1,#2.{\bigl\langle#1,#2\bigr\rangle}
\def\IP#1,#2.{\Bigl\langle#1,#2\Bigr\rangle}

\def\abs#1{\lvert#1\rvert}
\def\Abs#1{\bigl\lvert#1\bigr\rvert}
\def\ABs#1{\Bigl\lvert#1\Bigr\rvert}
\def\ABS#1{\biggl\lvert#1\biggr\rvert}

\newcommand{\za}{\ensuremath{\alpha}}
\newcommand{\zb}{\ensuremath{\beta}}
\newcommand{\zc}{\ensuremath{\psi}}

\newcommand{\zd}{\ensuremath{\delta}}

\newcommand{\ze}{\ensuremath{\epsilon}}
\newcommand{\zve}{\ensuremath{\varepsilon}}
\newcommand{\zvf}{\ensuremath{\varphi}}
\newcommand{\zf}{\ensuremath{\phi}}
\newcommand{\zF}{\ensuremath{\Phi}}

\newcommand{\zI}{\ensuremath{\infty}}
\newcommand{\zl}{\ensuremath{\lambda}}

\newcommand{\zm}{\ensuremath{\mu}}

\newcommand{\zt}{\ensuremath{\tau}}

\newcommand{\zx}{\ensuremath{\xi}}

\newcommand{\zs}{\ensuremath{\sigma}}

\def\z#1#2{\ifcase#1 {{\mathcal {#2}}}  %0=cal
\or {{\mathbf{#2}}}                    % 1=mathbold
\or { {\boldsymbol{#2}}}                  % 2=boldsymbol
\or{{\widetilde{#2}}}                   % 3=wide  tilde
\or {{\acute{#2}}}
\or {{\grave{#2}}}
\or {{\bar{#2}}}
\or {\dot{#2}}
\or {\overline{#2}}
\or {\underline{#2}}\fi}
 
\def\ZR{\ensuremath{\mathbb R}}
\def\ZZ{\ensuremath{\mathbb Z}}

\def\ZN{\ensuremath{\mathbb N}}

\def\ZS{\ensuremath{\mathbb S}}

\def\mid{\,:\,}

\def\md#1#2\emd{\ifx0#1
\begin{equation*} #2 \end{equation*}\fi  %  single line display, no number
\ifx1#1\begin{equation}#2\end{equation}\fi   % single line display, number
\ifx2#1\begin{align*}#2\end{align*}\fi   % aligned display, no number
\ifx3#1\begin{align}#2\end{align}\fi    % aligned display, number
\ifx4#1\begin{gather*}#2\end{gather*}\fi  % multiline, not align, no number
\ifx5#1\begin{gather}#2\end{gather}\fi   % multinline, not align
\ifx6#1\begin{multline*}#2\end{multline*}\fi  %  display too long for one line
\ifx7#1\begin{multline}#2\end{mutline}\fi  % as above, with numbers
}

\def\ind#1{ {\mathbf 1}_{#1}}
\def\size#1{\text{\rm eng($#1$)}}
\def\sh#1{\text{\rm sh}(#1)}

\def\exp#1 {{\text{\rm e}}^{#1}}  %   handling exponentials 
\def\T{ \ensuremath{{\mathbf T}}  }

\begin{document}

\title{ On the Hilbert Transform \\ and  Lacunary Directions in the Plane}
 \author{Michael T. Lacey\footnote{This work has been
  supported by an NSF grant, DMS--9706884.} 
\\ { Georgia Institute of Technology}  }
\date{}

\maketitle

\abstract{We show that the maximal operator below, defined initially for Schwartz
functions $f$ on the plane, extends to a bounded operator from  $L^p(\ZR^2)$ into itself for
$1<p<\zI.$
\md0
\sup_{k\in\ZZ} \ABs{ \text{p.v.} \int{}f(x-(1,2^k)y)\;\frac{dy}y}.
\emd
}

 \section{Introduction  }

For a smooth rapidly decreasing function $f$ on the plane and a nonzero vector
$v\in\ZR^2$, define
 \md0 H_vf(x):=\text{p.v}\int_{-\zI}^\zI f(x-yv)\frac{dy}y  ,
\emd
which is the one dimensional Hilbert transform of $f$ computed in the
direction $v$.  This definition is independent of the length of $v$.  
 For a countable collection of vectors $V\subset\ZR^2$, we define a maximal
function  \md0 H^Vf(x):=\sup_{v\in V}\abs{H_vf(x)}. \emd
  We shall prove

\bthm  t.V   Suppose $V=\{(1,a_k)\mid k\in\ZZ\}$ is  such that there is a $\zl>1$
so that for all $k$ \md0 0<a_{k+1}<a_k/\zl. \emd Then the operator $H^V$ extends
to a bounded operator from $L^p(\ZR^2)$ to itself for all $1<p<\zI$. \ethm

The theory of the Hilbert transform and maximal function of one variable are
closely intertwined.  And the maximal function version of the theorem above was
proved in a series of papers  \cite{m1,m2}, first in the $L^2$ case and last
of all for all $L^p$, $1<p<\zI$, \cite{m3}.  But the method of proof employed does
not seem to imply the theorem above.  

In a related matter, one can consider bounds on $H^V$ which depend only on the
cardinality of $V$.  On $L^2(\ZR^2)$, the bound of $\log \# V$ follows more or
less  immediately from the Rademacher-Menshov theorem.  But   maximal
function variants were only recently established by N.~Katz \cite{k1,k2}, using a
subtle range of ideas.

We shall prove the theorem above by invoking the BMO theory of the bidisk, as
developed by S.Y.~Chang and R.~Fefferman \cite{cf}.  This is conveniently done
via a combinatorial model of $H^V$, and once it is in place, a maximal
inequality can be established by way of an argument nearly devoid of the
geometry of the plane, as all the relevant geometric facts are already encoded
into the BMO theory.  

Indeed, the salient features of our argument are (1) a proper notion of
``energy" arising directly from a Bessel inequality, 
as used in e.g.~\cite{lt}, 
 (2) a closely related
notion of ``charge" and (3) a John--Nirenberg inequality.  Very little else is
needed to conclude a maximal inequality, making the proof adaptable to other situations with
the same attributes. This observation bears some resemblance to the method of
approach in N.~Katz' approach to the maximal function in arbitrary directions, 
\cite{k2}.

 \section{The Combinatorial Model}

Define the Fourier transform on $\ZR$ as $ \hat
f(\zx)=\int{}e^{-ix\zx}f(x)\;dx$.  We will use the same notation and a similar
definition  for the
Fourier transform on the plane $\ZR^2$.  Set $\ip f,g.:=\int f\bar g \, dx$. 
And by $A \seq B$, we mean that for some absolute constant $K$, $A\le{}K B$.

  We will replace the Hilbert
transform by $Pf(x)=\int_0^\zI e^{ix\zx}\hat f(\zx)\; d\zx$, which is Fourier
projection onto the positive frequencies of $f$.  $P$ is a linear combination
of the identity and the Hilbert transform.  Hence,  in complete analogy to the
definition of $H_v$,  we can define $P_vf$ as the one one-dimensional transform
applied  in the direction $v$ to the function $f$ defined on the plane. 

It suffices to consider maximal functions constructed from $P_v$, and more
particularly, it suffices to consider collections of vectors $V=\{(-1,a_k)\mid
k\in\ZN\} $ with $1/2<2^ka_{k}<1$ for $k\in\ZN$.    Throughout the rest of this
section we consider such collections $V$, and we define $P^V$ in complete
analogy to the definition of $H^V$.  

On the plane, view points as $x=(x_1,x_2)$, with the dual frequency variables
being $(\zx_1,\zx_2)$.  We define  
\md0 Bf(x)=\int_0^\zI\int_0^\zI \hat
f(\zx_1,\zx_2)e^{ix\cdot(\zx_1,\zx_2)}\;d\zx_1d\zx_2, 
\emd
 which is Fourier
projection onto $[0,\zI)^2$.  It suffices to consider the maximal function $
P^V Bf$.  Our purpose right now is to write $B$ as a limit of two different
combinatorial sums.  This will permit a corresponding decomposition of the
maximal operator.

For $j=1,2$, consider Schwartz functions $\zvf^j$ on $\ZR$ such that  \md4
\ind{[\frac78,\frac{13}8]}(\zx)\le\widehat{\zvf^1}(\zx)\le{}
\ind{[\frac34,\frac74]}(\zx), \\ \widehat{\zvf^2}(\zx)>0, \quad \zx>0, \\
\abs{\widehat{\zvf^2}(\zx)} {}\le{}  \min\{\abs\zx,\abs{\zx}^{-1}\}, \\
\zvf^2(x)\quad\text{is compactly supported.} \emd

Given a rectangle $R=r_1\times r_2$ in the plane, set  \md0
\zvf^j_R(x_1,x_2)=\prod_{k=1}^2\abs{r_k}^{-1/2}\zvf^k\Bigl(\frac{x_k-c(r_k)}{\abs{r_k}}\Bigr).
\emd Note that the $L^2$ norm of this function is independent of the choice of
$R$.

We shall consider classes of rectangles specified by
$\zl=(\zl_1,\zl_2)\in[1,2]^2$ and $y=\in\ZR^2$.  Set 
 \md0
\z0R_{\zl,y}=\Bigl\{y+\prod_{j=1}^2[\zl_jm_j2^{n_j},\zl_j (m_j+1)2^{n_j})\mid 
(m_1,m_2),(n_1,n_2)\in\ZZ^2\Bigr\}. \emd
We write $\z0R_{(1,1),(0,0)}$ as
$\z0R$, which is just (one choice of) all dyadic rectangles in the plane.

Notice that the functions $\{\zvf_R\mid R\in\z0R_{\zl,y}\}$ are obtained from
$\{\zvf_R\mid R\in\z0R\}$ by 
 dilating $x_1$ by a factor of $\zl_1$, $x_2$ by a
factor of $\zl_2$ with both dilations preserving the $L^2$ norm
 and then translating by $y$.
  Also note that there is no need to consider dilations by
factors greater than $2$, since a dyadic grid on $\ZR$ is invariant under
dilations by $2$.  

Define two operations by sums over these collections of rectangles.   \md0
C^j_{\zl,y}f(x):=\sum_{R\in\z0R_{\zl,y}}\ip f,\zvf^1_R.\zvf^j_R(x) ,\qquad j=1,2.
\emd We set
$C^j:=C^j_{(1,1),(0,0)}$.  

We use these to build two limiting representations of $B$, and to this this end
we note that $B$ is characterized, up to a constant multiple, as a non-zero
linear operator on $L^2(\ZR^2)$ that is (1) translation invariant (2) invariant
under dilations of both variables independently, and  (3)    $Bf=0$ if $\hat f$ is
not supported on $[0,\zI)^2$.  (This is suggested by the decomposition in
second section in \cite{lt}.)

 Define for $j=1,2$, 
\md0 B^jf(x)=\lim_{Y\to\zI}
\int\!\!\int_{D(Y)} C^j_{(2^{\zl_1},2^{\zl_2}),y}f\;\zm(d\zl_1,d\zl_2,dy), 
\emd 
where $\zm$ is normalized
Lebesgue measure on $D(Y):=[1,2)^2\times\{ y\mid \abs y<Y\}$. 
[Note that we are averaging with respect to multiplicative Haar measure in
 the dilation parameters $\zl_1$ and $\zl_2$.]
 The limit is seen to
exist for smooth compactly supported functions.  The operators $B^j$ extend to
bounded linear operators on $L^2$.  They are translation and  dilation invariant
since we average over all possible dilations and translations.  Moreover,
$B^jf=0$ if $\hat f $ is not supported on $[0,\zI)^2$.   It remains to show that
$B^j$, $j=1,2$ are non-zero operators and so are constant multiples of $B$.  Indeed, $B^1$
is easily seen to be positive semidefinite, by the choice of $\zvf^1$.  

For $B^2$,  observe that if we set $\zt_y\zvf(x)=\zvf(x-y)$ for functions on $\ZR$,
then we have  
\md2 \int_0^1\sum_{n=-\zI}^\zI\ip
f,\zt_{n+y}\zvf^1.\zt_{n+y}\zvf^2(x)\;dy
{}={}&\int_{-\zI}^\zI \int_{-\zI}^\zI
f(y)\overline{\zvf^1(-z)} \zvf^2(x-z)\; dy\,dz 
\\ {}={}&f*\zc(x) 
\emd  where
$\zc(x)=\int_{-\zI}^\zI \overline{\zvf^1(y)} \zvf^2(x+y)\; dy$.  Note that
$\hat\zc(\zx)=\overline{\hat{\zvf^1}(\zx)}\hat\zvf^2(\zx)$.  A direct
computation now shows that $B^2$ is non--zero.

We conclude that $C^j$ are constant multiples of $B$.  Hence, to prove our
theorem it suffices to prove a bound for  \md0 \sup_{v\in
V}\abs{C^2P_vC^1_{\zl,y}f},\qquad \zl\in[1,2)^2,\ y\in\ZR^2.
\emd We demonstrate a bound that is uniform in $\zl$ and $y$.  Then we can average these inequalities to conclude the same for
$\sup_{v\in V}\abs{B^2P_vB^1}$, which is sufficient.  

\medskip 
  Upon expansion of the term $C^2P_vC^1_{\zl,y}f$ we obtain the 
inner product $\ip  P_v\zvf_R^1,\zvf_{R'}^1.$ for $R\in\z0R$ and
$R'\in\z0R_{\zl,y}$.  This inner product 
is at most one in modulus. 
Moreover,  recall that we consider $v$ of the form $(-1,a)$, and write
$R=r_1\times r_2$.  Then 
\md0 P_{(-1,a)} \zvf_R^1={}  \begin{cases}
\zvf_R^1 &   \tfrac a2\abs{r_2}>\abs{r_1}   \\ 0 &  \abs{r_1}>2a\abs{r_2}
  \end{cases}  
\emd 
Thus, it is natural to  two cases, the
first being the sum taken over those rectangles with $P_v\zvf_R^1=\zvf_R^1$ and
the second is over those rectangles with $P_v\zvf_R^1\not\in\{0,\zvf_R^1\}$.  See
Figure~1.

% =========================================

\begin{center}
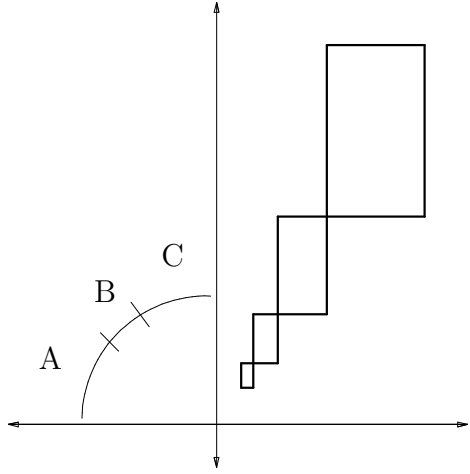
\begin{figure}
\begin{center}
\setlength{\unitlength}{0.00043333in}
\begingroup\makeatletter\ifx\SetFigFont\undefined%
\gdef\SetFigFont#1#2#3#4#5{%
  \reset@font\fontsize{#1}{#2pt}%
  \fontfamily{#3}\fontseries{#4}\fontshape{#5}%
  \selectfont}%
\fi\endgroup%
{\renewcommand{\dashlinestretch}{30}
\begin{picture}(5669,5759)(0,-10)
\drawline(2497.000,2122.000)(2421.394,2122.414)(2345.860,2119.090)
	(2270.583,2112.036)(2195.746,2101.270)(2121.533,2086.818)
	(2048.125,2068.716)(1975.702,2047.007)(1904.440,2021.744)
	(1834.514,1992.990)(1766.095,1960.814)(1699.350,1925.296)
	(1634.442,1886.522)(1571.531,1844.587)(1510.769,1799.594)
	(1452.305,1751.652)(1396.283,1700.878)(1342.839,1647.398)
	(1292.103,1591.341)(1244.201,1532.845)(1199.249,1472.052)
	(1157.357,1409.112)(1118.627,1344.178)(1083.154,1277.409)
	(1051.025,1208.968)(1022.319,1139.022)(997.105,1067.743)
	(975.445,995.305)(957.392,921.885)(942.990,847.662)
	(932.275,772.818)(925.273,697.536)(922.000,622.000)
\thicklines
\drawline(2872,997)(3022,997)(3022,1297)
	(2872,1297)(2872,997)
\drawline(3022,1297)(3322,1297)(3322,1897)
	(3022,1897)(3022,1297)
\drawline(3322,1897)(3922,1897)(3922,3097)
	(3322,3097)(3322,1897)
\drawline(3922,3097)(5122,3097)(5122,5197)
	(3922,5197)(3922,3097)
\blacken\thinlines
\drawline(142.000,577.000)(22.000,547.000)(142.000,517.000)(142.000,577.000)
\texture{55888888 88555555 5522a222 a2555555 55888888 88555555 552a2a2a 2a555555 
	55888888 88555555 55a222a2 22555555 55888888 88555555 552a2a2a 2a555555 
	55888888 88555555 5522a222 a2555555 55888888 88555555 552a2a2a 2a555555 
	55888888 88555555 55a222a2 22555555 55888888 88555555 552a2a2a 2a555555 }
\drawline(22,547)(5647,547)
\drawline(22,547)(5647,547)
\blacken\drawline(5527.000,517.000)(5647.000,547.000)(5527.000,577.000)(5527.000,517.000)
\drawline(1147,1672)(1372,1447)
\drawline(1522,2047)(1747,1747)
\thicklines
\blacken\thinlines
\drawline(2542.000,142.000)(2572.000,22.000)(2602.000,142.000)(2542.000,142.000)
\drawline(2572,22)(2572,5722)
\blacken\drawline(2602.000,5602.000)(2572.000,5722.000)(2542.000,5602.000)(2602.000,5602.000)
\put(397,1222){\makebox(0,0)[lb]{\smash{{{\SetFigFont{12}{14.4}{\rmdefault}{\mddefault}{\updefault}A}}}}}
\put(1072,2047){\makebox(0,0)[lb]{\smash{{{\SetFigFont{12}{14.4}{\rmdefault}{\mddefault}{\updefault}B}}}}}
\put(1897,2497){\makebox(0,0)[lb]{\smash{{{\SetFigFont{12}{14.4}{\rmdefault}{\mddefault}{\updefault}C}}}}}
\end{picture}
} 
\end{center}
\caption{\small  On the right are some rectangles $R$ in a collection $\z0R(k)$.  If a
vector $v$ is in arc $A$, then $P_v\zvf^1_R=0$, whereas if $v$ is in arc $C$, then
$P_v\zvf^1_R=\zvf^1_R$.  But if $v$ is in arc $B$, then $P_v\zvf^1_R$ need not be $0$ nor
$\zvf^1_R$. }\label{f.wR}
\end{figure}
\end{center}

% =========================================

\medskip 
 
 The second case concerns classes of rectangles which we define this way.  Recall
 that $V=\{(-1,a_k\mid k\in\ZN\}$ and set 
 \md0
 \z0R_{\zl,y}(k):=\{r_1\times r_2\in\z0R_{\zl,y}\mid
 2a_k\abs{r_2}\ge\abs{r_1}\ge\tfrac{a_k}2\abs{r_2}\}.
 \emd
 We again set $\z0R(k):=\z0R_{(1,1),(0,0)}(k)$.  These are the rectangles for
 which $P_{(-1,a_k)}\zvf^1_R$ need not be $0$ or $\zvf^1_R$.  We define 
 \md0
 \zF^j_{\zl,y,k}f:=\sum_{r\in\z0R_{\zl,y}}\ip  f,\zvf^1_R. \zvf^j_R,\qquad
 j=1,2, 
 \emd
 and set $\zF_k^j:=\zF^j_{(1,1),(0,0),k}$.  Then, in the case of $1<p<2$,
  the term to bound is 
 \md2
 \norm \sup_k\abs{ \zF^2_{\zl,y,k}P_{(-1,a_k)}\zF^1_k  f}.p.\le{}& 
 \NOrm\Bigl[\sum_k\abs{ \zF^2_{\zl,y,k}P_{(-1,a_k)}\zF^1_k  f}^2\Bigr]^{1/2}.p. 
 \\{}\seq{}&  \NOrm\Bigl[\sum_k\abs{  P_{(-1,a_k)}\zF^1_k  f}^2\Bigr]^{1/2}.p. 
 \\{}\seq{}&   \norm f.p.
 \emd
 The penultimate line is a vector valued Calderon-Zygmund inequality and the last
line follows from \l.P/ below. The case of $2\le p<\zI$ is even easier.

  \medskip

In the first case, there is the important point that the inner product $
\ip \zvf_R^1,\zvf_{R'}^1.$  will be zero unless the side lengths of
$R$ and $R'$ agree.  And if they do, the inner product will decay as a function
of the relative distance between $R\in\z0R$ and $R'\in\z0R_{\zl,y}$.  This relative distance will in 
addition be influenced by $y$, indeed, if $\abs y$ is considerably greater than
the side lengths of both rectangles, then it is the dominant term in
determining the relative distances between the two rectangles.

Thus, in seeking a useful quantitative estimate, it is useful to link the
translation parameter $y$ to the scales of the rectangles involved.  This we
shall explicitly do in this definition.  The maximal operator we control is
 \md1 \Label e.sum  
\sup_v\ABs{  \sum_{R\in\z0R}\ip f,\zvf^1_R
.\ip P_v\zvf_R^1,\zvf_{\zs(R)}^1.\zvf^2_{\zs(R)}}.
\emd
In this display, the map $\zs$ is given by 
\md1\Label e.zs 
\zs(r_1\times r_2)={}\zl_1 r_1\times \zl_2
r_2+(y_1\abs{r_1},y_2\abs{r_2})+(\zd_1({r_1}),\zd_2({r_2})),
\emd
in which $(\zl_1,\zl_2)\in(1,2]^2$ is fixed,  $y=(y_1,y_2)\in\ZR^2$ is fixed and
$\zd_j(\cdot)$, $j=1,2$ are two functions from the dyadic intervals on $\ZR$ into
$\ZR_+$, with $0\le\zd_j(r)\le{}\abs r$ for all dyadic intervals $r$.

Our purpose is to prove a bound on the $L^p$ norm of the operator in \e.sum/
which decays rapidly in $\abs y$.  This is possible because of the estimate 
 
\md0
\abs{\ip \zvf_R^1,\zvf_{\zs(R)}^1.}\le{}C(1+\abs{y_1}+\abs{y_2})^{-10}.
\emd
Therefore, the control of this part of the supremum will follow from the lemma
of the next section.   This completes our proof of the maximal theorem.

\def\sl#1{\text{\rm sl}(#1)}

\section{A Discrete Maximal Inequality}  For rectangles $R=r_1\times
r_2\in\z0R$ we set $\sl R:=\abs{r_1}/\abs{r_2}$.  

\bthm  l.max Let $y=(y_1,y_2)\in\ZR^2$ and let $\zs $ be as in \e.zs/.  
For all $1<p<\zI$ for arbitrary choices of signs $\{\zve_R\mid
R\in\z0R\}$, with $\zve_R\in\{\pm1\}$, the maximal operator below maps
$L^p(\ZR^2)$ into itself.  
 \md0 \sup_{s>0}\ABS{\sum_{ \substack{ R\in\z0R \\ 
\sl R\ge s}} \zve_R\ip f,\zvf_R^1.\zvf^2_{\zs(R)}} \emd 
The norm of the operator is at most $C_p(\abs{y_1}+\abs{y_2})^{3/p}$.  
\ethm 

\def\za#1{\langle f,\zvf_{#1}^1\rangle}

Let $\z0S$ denote an arbitrary finite subset of $\z0R$ and set  
\md0 
A^{\z0S}_{\max}={}  \sup_{a>0}\ABS{\sum_{ \substack{R\in\z0S\\\sl R>a}}\zve_R 
\za{R}
\zvf^2_{\zs(R)} }\emd
where $f\in L^p(\ZR^2)$ is a fixed function of norm one. By $A^{\z0S}$ denotes
the same sum over $R\in\z0S$, without the supremum over $a$. It
suffices to show that there is a constant $K_p$  independent of $f$ and 
  $\z0S\subset\z0R$, for which 
\md0 \abs{\{A^{\z0S}_{\max}>1\}}\seq(\abs{y_1}+\abs{y_2})^3. 
\emd
It follows that $A^{\z0S}_{\max}$ maps $L^p$ into  weak $L^p$ with norm at most
a constant times $(\abs{y_1}+\abs{y_2})^{3/p}$.  Interpolation then proves the
Lemma.

We define the shadow of $\z0S$ to be $\sh {\z0S}=\bigcup_{R\in\z0S}R$.  Since we
specified that $\zvf^2$ has compact support, it follows that 
 \md1 \Label e.contain  
\text{supp}(A^{\z0S})\subset{} \{M\ind{\sh{\z0S}}\ge\zd(1+\abs
{y_1}+\abs{y_2})^{-2}\} 
\emd 
for an absolute choice of 
 $\zd>0$, where $M$ is the strong maximal function. Thus,
$M$ can be defined as   \md0 Mg(x):=
\sup_{\substack{R\in\z0R\\x\in R}}\abs{2R}^{-1} \int_{2R} \abs{g(y)}\;dy. \emd

\def\zb#1{\frac{\ip f,\zvf_{#1}^1.}{\sqrt{\abs R}}}

Another definition we need is  the ``energy of a collection of rectangles $\z0S$"
\md0
\size{\z0S}:=\sup_{\z0S'\subset\z0S}\abs{\sh{\z0S'}}^{-1}
\NORm\Biggl[\sum_{R\in\z0S'}{\ABs{\zb R}}^2\ind R\Biggr]^{1/2}.1.  \emd
 This
quantity is related to the definition of BMO of the bidisk, as 
we have the equivalence $\size{\z0R}\simeq\norm f.BMO.$.  See \cite{cf,f}.
   We shall
specifically need the fact that 
\md0 \NORm\Biggl[\sum_{R\in\z0S}{\ABs{\zb
R}}^2\ind R\Biggr]^{1/2}.p.\seq \size{\z0S} \abs{\sh {\z0S}}^{1/p}. 
\emd
 This
is a manifestation of the John--Nirenberg  inequality for the BMO space.  See
\cite{cf}, or the concluding section of this paper.

Say that $\z0S $ has charge $\zd$ if $\size{\z0S}\le\zd$ and  
\md0
\sum_{R\in\z0S}\abs{\za R}^2\ge\frac {\zd^2} 4 \abs{\sh {\z0S}}. 
  \emd 
  There are 
two essential aspects of this definition.  For the first, if $\z0S$ has charge
$\zd$ then we may apply the Littlewood--Paley inequality, yielding    
\md0
 \zd \abs{\sh{\z0S}}^{1/p}\seq\NORm\Bigl[\sum_{R\in\z0S}{\ABs{\zb
R}}^2\ind R\Bigr]^{1/2}.p.  .
{}\seq\norm f.p. \emd

The second fact concerns a collection of rectangles    $\z0S$ of energy $\zd$.
It neccessarily contains a subset of charge $\zd$.  Suppose there are 
  two disjoint  subsets  $\z0S_1$ and $\z0S_2$ of $\z0S$,  of charge $\zd$. Then
$\z0S_1\cup\z0S_2$ also has charge $\zd$.   Therefore, $\z0S$ contains a (non--unique)
maximal subcollection of charge $\zd$.  

\medskip 

Having finished with definitions, the main argument begins. 
If we begin with a finite collection of rectangles $\z0S$, it has finite
energy.    We now reduce to the case in which $\size{\z0S}$ is at most one. 
Indeed, $\z0S$   is the union of collections $\z0S_j$ for $j\ge0$ with the energy
of $\z0S_0$ being at most one, and    the charge of $\z0S_j$ being $2^{j}$ for
all $j>0$.   As we are to prove a distributional inequality, it follows, after
an application of \e.contain/,  that we
need not consider those collections  $\z0S_j$ with  $j\ge0$.    

Thus we can assume that the energy of $\z0S$ is at most one. We shall show that
for any $q>p$,  
\md0 
\norm A^{\z0S}_{\max}.q.\le{}K_q {\size{\z0S}}^{1-2p/q} .
\emd 
That is, the estimate is now independent of $y_1$ and $y_2$ and depends on $p>1$
only through the implied constant in the inequality.  
This will conclude the proof of our lemma.

The construction which leads to this inequality begins now.  Define for
integers $j$  \md0 \z0S(j):=\{R\in\z0S\mid  \sl R\ge2^j\}. \emd Let $K_v$, for
integers $v\ge0$, be
subsets of $\ZZ$ such that 
 \begin{itemize}
 \item $K_0=\{j_0\}$ for some integer with $\z0S(j_0)=\emptyset$. 
  \item $K_v\subset K_{v+1}$ for all  $v$.
\item  For all $j\in\ZZ$ if $j_v\in K_v$ is the maximal element of $K_v$ less than or equal to
$j$ then    \md0 \size{\z0S(j)-\z0S(j_v)}\le2^{-v}. \emd \item  The cardinality
of $K_v$ is minimal, subject to the first two conditions. 
\end{itemize}

Notice that for each $j$ and $v$ with $j_v\not=j_{v-1}$, we have that
  $\size{\z0S(j_v)-\z0S(j_{v-1})}\ge2^{-v-1}$.
Then, we have for any integer $j$,  
\md0
 \Abs{A^{\z0S(j)}}\le{} \sum_{v=1}^\zI
\abs{A^{\z0S(j_v)}-A^{\z0S(j_{v-1})}} 
\emd And so it suffices to prove the
bound  
\md0 \Norm \sup_j \Abs{A^{\z0S(j_v)}-A^{\z0S(j_{v-1})}}.q.\seq{}
2^{-v(1-p/q)}, \qquad q>2p.
 \emd

For this last inequality, we again apply the decomposition of a set of
rectangles into subsets with charge.  Namely, for integers $w\ge v$ there are
collections $\ZS_w$ of subsets of $\z0S$ such that 
\begin{itemize} \item For all $w$, every  $\z0S'\in\ZS_w$ has charge $2^{-w}$.
\item For all $w$, the collections $\z0S'\in\ZS_w$ are pairwise disjoint.
\item  For each $v$ and $w\ge v$, there is an $\z0S'(v,w)\in\ZS_w$ such that  \md0
\z0S(j_v)-\z0S(j_{v-1})=\bigcup_{w=v}^\zI\z0S'(v,w). \emd 
\end{itemize} 
This is
achieved just by applying, to each collection of rectangles  
$\z0S(j_v)-\z0S(j_{v-1})$,
 the first decomposition above, adding in an additional step of pulling out 
 maximal subcollections of appropriate charge.

\medskip 

The collection $\ZS_w$ possesses the following property. 
\md0
\sum_{\z0S'\in\ZS_w}\abs{ \sh{\z0S'} }\seq{}2^{\max(2,p)w}.
\emd
Indeed, set 
 $\z0S^w:=\bigcup_{\z0S'\in\ZS_w}\z0S'$, which is a collection of 
rectangles of energy at most one, by the first step of our construction. 
In addition, set 
\md0
g:=\sum_{R\in\z0S^w}\ip f,\zvf_r^1.\zvf_r^2.
\emd
This function is supported on $\sh{\z0S^w}$.  And, 
as each $\z0S'\in\ZS_w$ has charge $2^{-w}$, 
\md0
\sum_{\z0S'\in\ZS_w}\abs{ \sh{\z0S'} }\seq{}2^{2w}
\lVert g\rVert_2^2.
\emd
And so we estimate the $L^2$ norm of $g$. 

At this point, we consider the case of $2\le{}p$.  In this case, 
we have in addition $\norm g.p.\seq\norm f.p.\seq1$, so that 
$\norm g.2.\seq\abs{\sh{\z0S^w}}^{\frac12-\frac1p}\norm g.p.$.  
And this proves our claim. 

Now, if $1<p<2$, then observe that $\norm g.BMO.\seq{}1$, so that 
\md0
\norm g.2.\seq{}\lVert g\rVert_p^{p/2}\lVert g\rVert_{BMO}^{1-p/2}\seq{}1.
\emd
And this finishes the proof of our observation. 

\medskip

To conclude,  we may estimate 
\md2 \bigl\lVert \sup_{\z0S'\in\ZS_w}
\abs{A^{\z0S'}}\big\rVert_q^q {}\le{}& \sum_{\z0S'\in\ZS_w} \bigl\lVert
{A^{\z0S'}}\big\rVert_q^q 
\\ {}\seq{}& 
 2^{-qw} \sum_{\z0S'\in\ZS_w} \abs{\sh{\z0S'}} 
\\ {}\seq{}& 
2^{-(q-2p)w} .
\emd
The proof of this lemma is done, as the $q$th root of the 
last estimate is summable in $w\ge{}v$, provided $q>2p$.

%%%%%%%%%%%%%%%%%%%%%%%%%%%%%%%   New section

\section{The Diagonal Terms}

We prove the lemma 

\bthm l.P  For $1<p<\zI$, we have the inequality 
\md0
\NOrm \Bigr[\sum_k \abs{ P_{(-1,a_k)}\zF^1_k f}^{p^*}\Bigr]^{1/p^*}.p.\seq\norm
f.p.,\qquad p^*:=\max(2,p).
\emd
\ethm 

The proof depends heavily on inequalities of  Littlewood--Paley type, namely 
\md0
\norm f.p.\simeq\NOrm \Bigl[\sum_k \abs{\zF^j_kf}^2\Bigr]^{1/2}.p.\simeq{}
\NOrm \Bigl[\sum_{R\in\z0R}\frac{\abs{\ip f,\zvf^j_R.}^2}{\abs R}\ind
R\Bigr]^{1/2}.p.,
\emd
where $1<p<\zI$ and $j=1,2$.   These are obtained from applications of the
ordinary Littlewood--Paley inequalities in each variable separately.  

\medskip 
The proof of the Lemma for $p\ge2$ is now at hand. 
\md2
  \sum_k \lVert P_{(-1,a_k)}\zF^1_k f\rVert_p^p
\seq{}&
\sum_k \lVert\zF^1_k f\rVert_p^p
\\{}\seq{}&\Bigl\lVert \Bigl[\sum_k \abs{\zF_k^1 f}^2\Bigr]^{1/2}\Bigr\rVert_p^p
\\{}\seq{}& \lVert f\rVert_p^p.
\emd

\medskip 

But for the case of $1<p<2$, we rely upon a more substantive approach.  We in
fact show that the operators 
\md0
T^*f:=\sum_{k=1}^\zI \zve_k  P_{(-1,a_k)}\zF^1_k f,\qquad \zve_k \in\{\pm 1\},
\emd
map $L^p$ into itself for $1<p<2$.  The constants involved are shown to be
independent of the choices of signs.  From this, the lemma follows.  

We use duality and prove the corresponding fact on the dual operator $T$.  One
readily checks that $T$ is bounded on $L^2$ and the point is to extend this fact
to $L^p$ for $p>2$.  Moreover, by the Littlewood--Paley inequalities, it
suffices to prove   the square function bound we state now. Set 
\md0
\zf_R:=P_{(-1,a_k)}\zvf^1_R, \qquad R\in\z0R(k),\ k\in\ZN.
\emd
Then we shall establish 
\md1 \Label e.sq  
\NOrm \Bigl[\sum_{R\in\z0R}\frac{\abs{\ip f,\zf_R.}^2}{\abs R}\ind
R\Bigr]^{1/2}.p.\seq\norm f.p.,\qquad 2<p<\zI.
\emd

This is clear for $p=2$ and so we see a second endpoint estimate with which to
interpolate.  The endpoint estimate is that the square function maps $L^\zI$
into BMO of the bidisk.  This fact requires that we prove the inequality 
\md1 \Label e.lbmo  
\sum_{R\subset U}\abs{\ip f,\zf_R.}^2\seq\abs U\lVert f\rVert_\zI^2,
\emd
for all open sets $U\subset \ZR^2$ and functions $f$. 

Once this is established, one can interpolate to deduce \e.sq/.  In fact the BMO
estimate directly supplies the restricted weak--type version of \e.sq/. [The
notions of energy and charge are relevant to this argument.  The details are left
to the reader.]  Then
standard interpolation supplies the inequality \e.sq/.

Our proof of \e.lbmo/ follows routine lines of argument, for the BMO theory of the
bidisk. We fix an open set $U$
and function $f$ bounded by $1$.  By the evident $L^2$ estimate, it suffices to
consider the case in which $f$ is supported off of the set $\{M\ind U>\frac12\}$.  
For a dyadic rectangle $R\subset U$ set 
\md0
\zm_R:=\sup\{\zm>0\mid \zm R\subset \{M\ind U>\tfrac12\}\}.
\emd
This quantity is at least $2$.  To a rectangle $R=r_1\times r_2$ we associate a set of
dyadic rectangles 
\md4
\z0S(R,j,\ell):=\{R'=r'_1\times r'_2\in\z0R\mid R'\subset R\ , \  r_j=r'_j 
\ ,\ \abs{R'}=2^{-\ell}\abs R\},\qquad j=1,2,\ \ell=1,2,\ldots
\emd
Set $\z0S(R)=\bigcup_{j=1}^2\bigcup_{\ell=1}^\zI\z0S(R,j,\ell)$. 
The principle fact to prove is that for each $R\subset U$, 
\md1\Label e.2bmo  
\sum_{R'\subset \z0S(R)}\abs{\ip \zf_{R'},f.}^2\seq\zm_R^{-1/20}\abs R.
\emd
This is so under the additional assumption that $f$ is bounded by $1$ and supported off of the
set $\{M\ind U>\frac12\}$.

That this proves \e.lbmo/ follows from an application of Journ\'e's Lemma, \cite{j}. The
details are left to the reader. [Virtually the only way to verify the Carleson measure
condition for a specific measure is through this basic lemma of Journ\'e.]

\bigskip

 There is a further reduction in \e.2bmo/ to make.  We assume that $\zm>2$, that  $f$ 
 is bounded by $1$ and supported on $2\zm R-\zm R$.  Then we show that 
 \md1\Label e.22bb  
 \sum_{R'\subset \z0S(R)}\abs{\ip \zf_{R'},f.}^2\seq\zm^{-1/8}\abs R.
\emd
 This proves our desired inequality, as is easy to see.
 This inequality  depends upon specific properties of the
transformation $P_v$.

 Now, for $R'\in\z0R(k)$, the function
$\zf_{R'}=P_{(-1,a_k)}\zvf^2_{R'}$ need not be $0$ or $\zvf^2_{R'}$, and we have the
following estimate. 
\md1\Label e.ptwise  
\abs{\zf_{R'}(x)}\seq\abs{R'}^{-1/2}\{1+\abs{r'_1}[(x-c(R'))\cdot (-1,a_k)]+{}
(\abs{r'_2}[(x-c(R'))\cdot (a_k,1)])^{200}\}^{-1}
\emd
A standard calculation verifies this, recalling that the kernel associated to $P$---and
hence $P_v$---has decay of order $1/y$.  It then follows that if in addition
$R'\in\z0S(R,j,\ell)$, we then have 
\md0
\int_{2\zm R-\zm R}\abs{\zf_{R'}}\; dx\seq\sqrt{\abs{R'}}\frac{\log \zm}\zm.
\emd
We use this estimate for all $R'\in\z0S(R,j,\ell)$, with $j=1,2$ and $\ell\le\sqrt \zm$.
 As there are $2^\ell$ member of $\z0S(R,j,\ell)$, we see that 
 \md0
 \sum_{j=1}^2\sum_{\ell=1}^{\sqrt\zm}\sum_{R'\subset \z0S(R,j,\ell)}\abs{\ip 
 \zf_{R'},f.}^2\seq{} 
 \frac{\abs{R}}{\sqrt{\zm}}.
 \emd
 
 \smallskip 
 
 It remains to consider the case of $\ell\ge\sqrt\zm$.  Let us consider $j=1$, the case
 of $j=2$ being simmilar.  Without loss of generality, we can assume that $R$ is
 centered at the origin.  The difficulty arises from the poor decay of the functions $\zf_{R'}$
for $R'\in\z0R(k)$ in the direction $(-1,a_k)$.  However, these directions are now
localized in the direction $(-1,0)$ since $\ell\ge\sqrt\zm$. 
[By assumption $a_k\simeq2^{-k}$ and $\ell\ge{}k$.]
 Thus we break up the
support of $f$ in this way. 
\md4
V_1:=\{(-2\zm\abs{r_1},-\zm\abs{r_1})\cup(\zm\abs{r_1},2\zm\abs{r_1})\}\times(-\sqrt\zm
\abs{r_2},\sqrt\zm \abs{r_2}),
\\
V_2:=2\zm R-\zm R-V_1.
\emd
The critical inequality is that if $g$ is supported on $V_1$, then for all
$\ell\ge\sqrt\zm$, 
\md1\Label e.22 
\sum_{R'\in\z0S(R,1,\ell)}\abs{\ip \zf_{R'},g.}^2\seq\zm^{-9/10}\lVert g\rVert_2^2.
\emd
Indeed, from \e.ptwise/, we see the inequalities 
\md4
\norm \zf_{R'}.L^2(V_1).\seq\zm^{-1},\qquad R'\in\z0S(R,1,\ell),
\\
\int_{V_1}\abs{\zf_{R'}\zf_{R''}}\;dx\seq\Biggl[\frac{\text{dist}(R',R'')}
{2^{-\ell}\abs{r_2}}
\Biggr]^{-100},\qquad R',R''\in\z0S(R,1,\ell).
\emd
The estimate \e.22/ is a direct consequence of these two observations.  
Applying this to a bounded function $f$ supported on $V_1$, we see that $\norm
f.2.\le{}
\norm f.\zI.\sqrt{\abs{V_1}}\seq\zm^{3/4}\abs{R}$.  This then is consistent with 
 our claim \e.22bb/.
 
 The remaining case concerns bounded functions $f$ supported on $V_2$.  But in this
 case, the decay of the functions $\zf_{R'}$ is much better.  From \e.ptwise/ and the
 definition of $V_2$ we see that 
 \md0
 \int\abs{\zf_{R'}}\; dx\seq\sqrt{R'}2^{-10\ell},\qquad R'\in\z0S(R,1,\ell).
 \emd
 This is more than enough to conclude \e.22bb/ for bounded functions $f$ supported on $V_2$,
 finishing the proof of that inequality.

\appendix 

\section{Appendix: Carleson Measures}

We include a proof of some results related to the BMO theory of the bidisk
 cited above.  Our results will be slightly more general than we need. Let us 
  begin with the
John--Nirenberg Lemma.  For the collection of dyadic rectangles $\z0R$ in the
plane, let $a:\z0R\to\ZR_+$ be a map.  Define 
\md0
\norm a.CM,p.:=\sup_{U\subset\ZR^2}\abs{U}^{-1}\NOrm \sum_{R\subset U}a(R)\ind R
.p.,\qquad 0\le{}p<\zI,
\emd
where the supremum is over all open sets $U\subset\ZR^2$. 
This is a possible definition of the norm of a Carleson measure.
The John--Nirenberg inequality asserts that all of these possible definitions
are equivalent, up to constants. 

\bthm l.jn  For all $0\le{}p,q<\zI$, we have 
\md0
\norm a.CM,p.\seq \norm a.CM,q.
\emd \ethm

\begin{proof}
For open sets $U\subset\ZR^2$ define 
$$
F_U(x):=\sum_{R\subset U}a(R)\ind R(x).
$$

It suffices to prove the inequality above for a restricted range of $p$ and $q$.
We begin with the case of $\norm a.CM,1.\seq\norm a.CM,p.$, for some $0<p<2$. 
It suffices to fix a choice of 
 $a$ with $\norm a.CM,p.\le1$, and  $\text{supp}(F_{\z0R})=U$ has finite mesure.
 We need only show that 
 \md0
 \int F_U\; dx\seq\abs{U}
 \emd
 And to this end, it suffices to demonstrate that there is an open set
 $V\subset U$, with $\abs V<\abs{U}/2$, for which 
 \md0
 \int F_U\seq\abs{U}+\int F_v\; dx.
 \emd
 It is clear that this inequality can then be inductively applied to $V$ to
 yield the proof of the desired inequality.

We define $V$ as follows. For some constant $0<\ze<1/2$, set $E:=\{
F_U>\ze^{-2/p}\}$,
and $
V:=\{ M\ind E>\ze \}$.  By the boundedness of the strong maximal
function, for $\ze$ appropriately small, the measure of $V$ is at most one--half
the measure of $U$.  

But at the same time, if $R$ is a dyadic rectangle with $R\not\subset V$, then
$\abs{R\cap E}<\abs R/2$.  Hence, 
\md2
\int (F_U-F_V)\; dx={}&\sum_{\substack{R\subset U\\R\not\subset V}}a(R)\abs
R
\\{}\seq{}& \sum_{\substack{R\subset U\\R\not\subset V}}a(R)\abs{
R\cap E^c}
\\{}\seq{}& \int_{E^c}(F_U)\; dx
\\{}\seq{}& \int F_U^p\; dx
\\{}\seq{}& \abs U.
\emd
This follows since we have an upper bound of $F_U$ off of the set $E$. This case
has been proved.

\medskip 

We now turn to the estimate $\norm a.CM,p.\seq\norm a.CM,1.$,  for choices of
$1<p<\zI$.  This is all that remains to be done.  Indeed this is the
case that is explicitly proved in \cite{cf}, but we include the details for the
convenience of the reader.

Due to the recursive nature of the definition of energy, it suffices to prove
the following. Fix a choice of 
 $a$ with $\norm a.CM,p.\le1$, and  $\text{supp}(F_{\z0R})=U$ has finite mesure.
Then, there is an open set $V\subset U$, with $\abs V<\abs{U}/4$, for which
\md0
\norm F_U .p.\seq\abs{U}^{1/p}+\norm G_V.p..
\emd

This is done by way of duality.  Thus  let $p'$ be the conjugate index to $p$ 
and
select  a non--negative $h\in L^{p'}$  of norm one so that $\ip F_U,h.=\norm 
F_U.p.$.  The open set $V$ is then 
 \md0
 V:=\bigcup_{R\in\z0R}\bigl\{R\mid \frac1{\abs{R}}\int_{R}h\;dx>\zl\bigr\},
 \emd
 where we choose $\zl>0$ momentarily.
 
 By the boundedness of the strong maximal function,  
 \md0
 \abs V\le{}C_r\zl^{-p'}{\lVert h\rVert_{p'}^{p'}}=\tfrac12\abs U
 \emd
 if we take $\zl\simeq\abs{U}^{1/p'}$.  But then 
 \md2
 \norm F_U.p.=\ip F_U,h.\le{}&
                \norm F_V.p.+\sum_{\substack{R\in\z0R\\ R\not\subset V}}
                  {\abs{a_s}}   \int_{R}h\;dx
  \\
  {}\le{}&\norm F_V.p.+\zl\sum_{s\in\T} \abs{a_s} 
  \\
  {}\le{}&\norm F_V.p.+\zl\abs U,
  \emd
  which proves our inequality by the  choice of $\zl$.
 
\end{proof}

The second topic is that of Journ\'e's Lemma, which in any of it's various forms
must be stated in terms of these quantities.  Fix an open set $U$ and a
rectangle $R\subset U$.  Then 
\md0
\zm_R:=\sup \{\zm \mid \zm R\subset \{M\ind U>1/2\}\}.
\emd 
 Specifically, in this paper we
assumed this lemma. [For a more precise result, see \cite{j}.]

\bthm l.journe Fix $\ze>0$. Let $a\mid \z0R\to\ZR_+$ be such that for all open sets
$U$ and all dyadic $R\subset U$, 
\md0
\sum_{R'\subset R}a_{R'}\seq\zm^{-\ze}\abs R.
\emd
Then $\norm a.CM,1.\seq1$.  
\ethm

\begin{proof}    Notice that if $\z0R'$ is a
collection of rectangles for which 
\md0
\abs{R\cap R'}<\tfrac 12(\abs{R}\wedge\abs{R'}), \qquad R,R'\in\z0R'
\emd
then 
\md0
\sum_{R\in\z0R'}\abs{R}\le2\abs{\cup\{R\mid R\in\z0R'\}}.
\emd
Our obective is to arrange the collection of rectangles into subcollections
which are ``nearly disjoint" in this sense, and for which $\zm_R$ is
approximtely the same energy.

   For integers $k\ge0$, let $\z0R_k$ be those dyadic rectangles
which satisfy (1) $R\subset U$, (2) $R$ is maximal among all   rectangles
satisfying (1), (3)
and $2^k\le\zm_R<2^{k+1}$.    Then let $\z0R'_k$ be a subcollection of $\z0R$ in which the
lengths of the two sides of $R$ are restricted to be in $2(k+1)\ZZ+j$ for the first
side and $2(k+1)\ZZ+j'$, with $0\le{}j,j'<2(k+1)$.  Certainly there are at most
$4(k+1)^2$ such subcollections $\z0R_k'$.  But, by maximality,
 these collections are ``nearly disjoint" in the sense of the previous
 paragraph.  Hence, 
 \md0
 \sum_{R\in\z0R'_k}\sum_{R'\subset R}a_{R'}
 {}\le{} 2^{-\ze k}\sum_{R\in\z0R'_k}\abs R\le{} 2^{-\ze k}\abs{U}.
 \emd
 This estimate is summable over the $4(k+1)^2$ possible choices of $\z0R'_k$ and
 over $k\ge0$.
 \end{proof}

\bigskip
 {\parindent=0pt\baselineskip=12pt\obeylines
Michael T. Lacey 
School of Mathematics
Georgia Institute of Technology
Atlanta GA 30332
\smallskip
\tt lacey@math.gatech.edu
\tt http://www.math.gatech.edu/\~{}lacey }
 
\end{document}